\documentclass[a4paper,myheadings]{amsart}

\usepackage{amssymb}
\usepackage{amsthm}
\usepackage{amsmath}
\usepackage{eucal}
\usepackage{mathrsfs}

\theoremstyle{definition}

\newcommand{\cal}{\mathcal}

\newcommand{\Con}{\ensuremath{\mathrm{C}}}
\newcommand{\Cinf}{\ensuremath{\mathrm{C}^\infty}}
\newcommand{\D}{\ensuremath{{\cal D}}}
\renewcommand{\S}{\ensuremath{{\cal S}}}
\newcommand{\E}{\ensuremath{{\cal E}}}

\newcommand{\mb}[1]{\ensuremath{\mathbb{#1}}}
\newcommand{\N}{\mb{N}}

\newcommand{\R}{\mb{R}}

\newcommand{\cl}[1]{\ensuremath{[#1]}}

\newcommand{\G}{\ensuremath{{\cal G}}}

\newcommand{\NN}{\ensuremath{{\cal N}}}
\newcommand{\Ginf}{\ensuremath{\G^\infty}}

\renewcommand{\d}{\ensuremath{\partial}}
\newcommand{\diff}[1]{\frac{d}{d#1}}

\newfont{\bl}{msbm10 scaled \magstep2}

\newtheorem{thm}{Theorem}[section]

\newtheorem{prop}[thm]{Proposition}

\newtheorem{cor}[thm]{Corollary}
\newtheorem{rem}[thm]{Remark}

\newcommand{\beq}{\begin{equation}}
\newcommand{\eeq}{\end{equation}}

\newcommand{\map}{\ensuremath{\rightarrow}}

\newcommand{\FT}[1]{\widehat{#1}}

\newcommand{\F}{\ensuremath{{\cal F}}}

\newcommand{\inp}[2]{\langle #1 | #2 \rangle}  
\newcommand{\notmid}{\mid\kern-0.5em\not\kern0.5em}

\newcommand{\norm}[2]{{\| #1 \|}_{#2}}

\newcommand{\ltwo}[1]{\norm{#1}{L^2}}

\newcommand{\linf}[1]{\norm{#1}{L^\infty}}

\newcommand{\al}{\alpha}
\newcommand{\be}{\beta}
\newcommand{\ga}{\gamma}

\newcommand{\de}{\delta}
\newcommand{\eps}{\varepsilon}

\newcommand{\la}{\lambda}

\newcommand{\om}{\omega}

\newcommand{\inv}{^{-1}}

\newcommand{\floor}[1]{\ensuremath{\lfloor #1 \rfloor}}
\renewcommand{\Re}{\ensuremath{\text{Re}}}
\renewcommand{\Im}{\ensuremath{\text{Im}}}

\newcommand{\ovl}[1]{\overline{#1}}

\newcommand{\adjoint}[1]{\ensuremath{{#1}^{\textstyle{*}}}}

\newcommand{\Gtwo}{\ensuremath{\G_{L^2}}}

\newcommand{\wt}[1]{\ensuremath{\widetilde{#1}}}

\newcommand{\dq}{d\hspace{-0.4em}{ }^-\hspace{-0.2em}}
\renewcommand{\S}{\mathscr{S}}

\begin{document}

\title{First-order hyperbolic pseudodifferential equations with generalized symbols}

\author{G\"{u}nther H\"{o}rmann}\thanks{Supported by FWF grant P14576-MAT}
\address{Institut f\"ur Mathematik, Universit\"at Wien}
\curraddr{Institut f\"ur Technische Mathematik, Geometrie und
Bauinformatik, Universit\"at Innsbruck}
\email{guenther.hoermann@univie.ac.at}

\date{\today}


\maketitle \markboth{G\"unther H\"ormann}{Hyperbolic pseudodifferential
equations with generalized symbols}

\begin{abstract}
We consider the Cauchy problem for a hyperbolic pseudodifferential operator
whose symbol is generalized, resembling a representative of a Colombeau
generalized function. Such equations arise, for example, after a
reduction-decoupling of second-order model systems of differential equations
in seismology. We prove existence of a unique generalized solution under
log-type growth conditions on the symbol, thereby extending known results for
the case of differential operators (\cite{LO:91, O:89}).
\end{abstract}

AMS 2001 Mathematical Subject Classes: 46F30, 35S10 

\section{Introduction}

This paper establishes existence and uniqueness of a generalized solution to
the scalar hyperbolic pseudodifferential Cauchy problem
\begin{align}
    \d_t u + A(t,x,D_x) u &=  f  \quad \text{when $t\in (0,T)$}, \label{PsDE}\\
    u(0) &= g. \label{initialvalue}
\end{align}
The data $f$ and $g$ are Colombeau generalized functions and $A$ is a
generalized pseudodifferential operator of order $1$. Its symbol is represented
by a family of smooth regularizations, which may (but need not) be convergent
to a distributional symbol. Problem (\ref{PsDE}-\ref{initialvalue}) represents
an extension of the scalar case of the partial differential equations
considered by Lafon and Oberguggenberger in \cite{LO:91,O:89}.

One may think of problem (\ref{PsDE}-\ref{initialvalue}) as resulting from a
system of second-order (partial differential) equations by reduction to
first-order followed by a decoupling into scalar equations (cf.\ \cite[Section
IX.1]{Taylor:81}). This is a standard technique in applications, for example,
in mathematical geophysics, where one decouples the modes of seismic
propagation and subjects these to further refined analysis (cf.\
\cite{SdH:02}). As they stand, these reduction-decoupling methods are
rigorously applicable in the case of models with smooth coefficients or
symbols, but cease to be well-defined under the realistic assumptions of only
measurable (bounded) coefficients, which are to represent the elastic or
acoustic properties of the earth's subsurface. Moreover, the initial value and
the right-hand side are distributions corresponding to the original seismic
source and force terms, which are, by nature, strongly singular, e.g.,
delta-like. If the original model coefficients are replaced by regularizations,
then we may carry out all transformations within algebras of generalized
functions from the outset and arrive at (\ref{PsDE}-\ref{initialvalue}) in a
well-defined way. The purpose of the current paper is to investigate the
general feasibility of rigorously solving the resulting decoupled, so-called
one-way wave equation, by generalized functions. Future work will be devoted to
the regularity analysis of the solutions and their asymptotic relations with
distributions.

A word on \emph{conventions and notations concerning the Fourier transform:} if
$u$ is a temperate distribution on $\R^n$ we denote its Fourier transform by
$\FT{u}$ or $\F u$; occasionally, when several variables and parameters are
involved, we write expressions of the form $\F_{x\to\xi} (u(y,x))$ to indicate
that the transform acts on the partial function (or distribution) $u(y,.)$ and
yields a function (or distribution) in $(y,\xi)$; the integral formulas for the
transforms follow the convention $\F u (-x)/(2\pi)^n =  
\F\inv u(x) = \int \exp(i x\xi)\, u(\xi)\, \dq\xi$, where 
$\dq\xi = d\xi/ (2\pi)^n$. 

Subsections 1.1-3 serve to review Colombeau theory, fix our notations for
generalized symbols, and also recall the corresponding result on the Cauchy
problem for hyperbolic differential operators with generalized coefficients.
Section 2 establishes precise energy estimates, which are at the heart of the
existence and uniqueness proof for the Cauchy problem presented in Section 3.
Finally, under additional assumptions on the symbol and data regularity, we
are able to draw some conclusions about the solution regularity which are revealed by the technique of the existence proof itself.

\subsection{Colombeau algebras of generalized functions}

We will set up and solve the problem in the framework of algebras of
generalized functions introduced by Colombeau in
\cite{Colombeau:84,Colombeau:85}. More specifically, we will work in a variant
which is based on	 
$L^2$-norm estimates as introduced in \cite{BO:92}. We will recall the
definition and basic properties below. As general references 
and for discussions of the overall properties of Colombeau
algebras we refer to the literature (e.g. \cite{Colombeau:85,GKOS:01,O:92}).

We consider the space-time domain $X_T := \R^n\times (0,T)$. The basic objects
defining our generalized functions are regularizing families $(u_\eps)_{\eps
\in (0,1]}$ of smooth functions $u_\eps \in H^\infty(X_T)$ for $0 < \eps \leq
1$, where $H^\infty$ denotes the intersection over all Sobolev spaces. To
simplify the notation, we shall write $(u_\eps)_\eps$ in place of
$(u_\eps)_{\eps \in (0,1]}$ throughout. We single out the following
subalgebras:
\vspace{2mm}\\
{\em Moderate families}, denoted by $\E_{M,L^2}(X_T)$, are defined by the
property:
\[
    \forall \alpha \in \N_0^n\, \exists p \geq 0:\;
    \ltwo{\d^\alpha u_\eps} = O(\eps^{-p})\ \rm{as}\ \eps \to 0\,.
\]
{\em Null families}, denoted by $\NN_{L^2}(X_T)$, are the families in
$\E_{M,L^2}(X_T)$ having the following additional property:
\[
    \forall q \geq 0:\; \ltwo{ u_\eps} = O(\eps^q)\ \rm{as}\ \eps \to 0\,.
\]
Hence moderate families satisfy $L^2$-estimates with at most polynomial
divergence as $\eps \to 0$, together with all derivatives, while null families
vanish faster than any power of $\eps$ in the $L^2$-norms. For the latter,
one can show that, equivalently, all derivatives
satisfy estimates of the same kind (cf.\ \cite{Garetto:02b}). The null families
form a differential ideal in the collection of moderate families. The {\em
Colombeau algebra} $\Gtwo(X_T)$ is the factor algebra
\[
   \Gtwo(X_T) = \E_{M,L^2}(X_T) / \NN_{L^2}(X_T)  \,.
\]
(The notation in \cite{BO:92} is $\G_{2,2}$, and correspondingly for moderate
and negligible nets, where the variability of $L^q$-norms in the definitions
was essential.)
The algebra $\Gtwo(\R^n)$ is defined in exactly the same way and its elements
can be considered as elements of $\Gtwo(X_T)$. On the other
hand, as explained in \cite[Remark 2.2(i) and Definition 2.8]{BO:92}, the
restriction of a generalized function from $\Gtwo(X_T)$ to $t=0$ is
well-defined: for any representative $(u_\eps)_\eps$  in $\E_{M,L^2}(X_T)$ we
have $u_\eps \in \Cinf(\R^n\times [0,T])$ (i.e., smoothness up to the boundary
of the time interval) and that $(u_\eps(.,0))_\eps$ belongs to
$\E_{M,L^2}(\R^n)$. We use the bracket notation $\cl{\ .\ }$ to denote the
equivalence class in $\Gtwo$.

Distributions in $H^{-\infty}(\R^n) = \bigcup_{s\in\R} H^s(\R^n)$ are embedded
in $\Gtwo(\R^n)$ by convolution:  $\iota(w) = \cl{(w * (\rho_\eps))_\eps}$,
where
\begin{equation}
   \rho_\eps(x) = \eps^{-n}\rho\left(x/\eps\right)
               \label{molli}
\end{equation}
is obtained by scaling the fixed mollifier $\rho$, i.e., a test function $\rho
\in \S(\R^n)$ of integral one with all moments (of order $1$ and higher)
vanishing. This embedding renders $H^\infty(\R^n)$ a faithful subalgebra (cf.\
\cite[Theorem 2.7(ii)]{BO:92}). In fact, given $f \in H^\infty(\R^n)$, one can
define the corresponding element of $\Gtwo(\R^n)$ by $\cl{(f)_\eps}$ (with
representative independent of $\eps$). In the same way we may consider
$H^\infty(X_T)$ a faithful subalgebra of $\Gtwo(X_T)$. 

Some Colombeau generalized functions behave macroscopically like a
distribution. 
We say that $u = \cl{(u_\eps)_\eps}\in\Gtwo$ is \emph{associated with the
distribution} $w\in\D'$, denoted by $u \approx w$, if $u_\eps \to w$ in $\D'$
as $\eps \to 0$. 

Intrinsic regularity theory for Colombeau generalized functions has been
a subject of active research. Its foundation is \cite[Section
25]{O:92} with the definition of the subalgebra $\Ginf$ of $\G$, which plays
the same role for $\G$ as $\Cinf$ does within $\D'$. The basic idea is to
couple the generalized regularity notion to uniform $\eps$-growth in all
derivatives and it leads to the important compatibility relation
 \[
	\Ginf \cap \D' = \Cinf.
 \]
Similarly, we define here the subalgebra $\Gtwo^\infty$ of \emph{regular
elements} of $\Gtwo$ by the following condition: $u = \cl{(u_\eps)_\eps}
\in\Gtwo$ belongs to $\Gtwo^\infty$ if and only if 
 \beq\label{regularity}
   \exists p \geq 0 \, \forall \alpha \in \N_0^n:\;
    \ltwo{\d^\alpha u_\eps} = O(\eps^{-p})\ \rm{as}\ \eps \to 0\,.	
 \eeq
Observe that $p$ can be chosen uniformly over all $\al$. In particular, if
$u_\eps = v * \rho_\eps$ with $v\in H^\infty$, then $p = 0$ is possible when
we let fall all derivatives on the factor $v$.

Concerning sources for recent and related research in Colombeau theory, with a
diversity of directions, including such topics as pseudodifferential 
operators with generalized symbols, regularity theory, and
microlocal analysis of nonlinear singularity propagation we refer to
\cite{Garetto:02,HdH:01,HK:01,HO:02,HOP:03,NPS:98}.

\subsection{Generalized pseudodifferential operators}

For comprehensive theories of approaches to pseudodifferential operators with
Colombeau generalized functions as symbols we may refer to the recent
literature on the subject \cite{Garetto:02,GGO:03,NPS:98}. However, the purpose
of the present paper is to present a short and self-contained discussion of the
solution to the hyperbolic pseudodifferential Cauchy problem. Therefore we do
not need to call on the full theory of generalized symbol classes, mapping
properties, and symbol calculus, as it has been extended systematically and
with strong results in \cite{Garetto:02,GGO:03}. Nevertheless, this background
will be substantial in further development, refinements, and applications of
the current work, in particular, concerning regularity theory and microlocal
analysis.

We will use families of smooth symbols satisfying uniform estimates with
respect to the $x$ (and $t$) variable as described in
\cite{Hoermander:V3,Kumano-go:81}.  To fix notation, let us review the
definition. A complex valued function $a\in\Cinf(\R^n\times\R^n)$ belongs to
the symbol class $S^m$ of order $m\in\R$ if for all $(\al,\be)\in\N_0^{2n}$
\beq\label{symbdef}
    c^m_{\al,\be}(a) := \sup\limits_{(x,\xi)\in\R^{n+p}}
        (1 + |\xi|)^{-m + |\al|} |\d_\xi^\al \d_x^\be a(x,\xi)| < \infty.
\eeq
$S^m$ is a Fr\'echet space when equipped with the semi-norms
\beq
    q^m_{k,l}(a) := \max\limits_{|\al| \leq k, |\be| \leq l} c^m_{\al,\be}(a),
\eeq
a notation we will make use of freely in several estimates in the sequel.
(Observe that compared to the semi-norms and notation used in
\cite{Kumano-go:81} we have
$|a|^{(m)}_l = \max \{ q^m_{k,r}(a) : k+r \leq l \} $.)
In fact, we will use symbols which depend smoothly on time, considered as a
parameter. More precisely, we consider the space of symbols $a(t,x,\xi)$ where
$a\in\Cinf([0,T],S^m)$ (i.e., each $t$-derivative on $(0,T)$ is continuous
into $S^m$ up to the boundary $t=0$ and $t=T$) with the semi-norms
\beq
    Q^m_{j,k,l}(a) := \max\limits_{0\leq i \leq j}
        \sup\limits_{t\in[0,T]} q^m_{k,l}(\d_t^i a(t,.,.)).
\eeq

By a \emph{generalized symbol} we mean a family
$(a_\eps)_{\eps\in [0,1)}$ of smooth symbols in $S^m$ (the
same $m$ for all $\eps$) which satisfy moderate semi-norm estimates, i.e., for
all $k$ and $l$ in $\N_0$ there is $N\in\N_0$ such that
 \beq\label{modsymb}
    q^m_{k,l}(a_\eps) = O(\eps^{-N}) \qquad (\eps\to 0).
 \eeq Generalized symbols with parameter $t\in[0,T]$ are given by families
$(t,x,\xi) \mapsto  a_\eps(t,x,\xi)$ ($\eps\in(0,1]$) such that
$a_\eps\in\Cinf([0,T],S^m)$ with moderate semi-norm estimates: for all $j$,
$k$, $l$ in $\N_0$ there is $N\in\N_0$ such that
 \beq\label{modsymb_t}
    Q^m_{j,k,l}(a_\eps) = O(\eps^{-N}) \qquad (\eps\to 0).
 \eeq
Obviously, no major changes would be required to incorporate more general types
of symbols, especially the H\"ormander's classes $S^m_{\rho,\de}$ would mainly
require changes in notation (at least when $0 \leq \rho < \de < 1$).

Let $(a_\eps)_\eps$ be a generalized symbol with parameter $t\in[0,T]$. We
define the corresponding linear operator
 \[
    A : \Gtwo(X_T) \map \Gtwo(X_T)
 \]
 in the following way. On the representative level, $A$ acts as the diagonal operator
 \[
    (u_\eps)_\eps \mapsto (a_\eps(t,x,D_x) u_\eps)_\eps \qquad
        \forall (u_\eps)_\eps\in \E_{M,L^2}(X_T).
 \]
Here, $a_\eps(t,x,D_x)$ acts as an operator in the $x$ variable with
parameter $t$. The moderateness of $(a_\eps(t,x,D_x)
u_\eps)_\eps$ follows from (\ref{modsymb}) and the fact that
operator norms of $\d_t^i \d_x^\be \circ a_\eps(x,t,D_x)$ on Sobolev spaces
are bounded (linearly) by finitely many semi-norms of the symbol
(cf.\ \cite[Ch.\ 3, Theorem 2.7]{Kumano-go:81}). In the same way, it
follows that null families are mapped into null families, so that $A$ is
well-defined on equivalence classes. We call $A$ the \emph{generalized
pseudodifferential operator} with generalized symbol $(a_\eps)_\eps$.

\subsection{Review of hyperbolic partial differential equations with
generalized coefficients}

We briefly review the situation for symmetric hyperbolic systems of partial
differential equations in Colombeau algebras. The heart of this theory was
developed in \cite{LO:91,O:89}, from where we recall the main result on the
Cauchy problem.

The theory is placed in $\G$ instead of $\Gtwo$, i.e., the data $f$, $g$ as
well as all coefficients satisfy asymptotic local $L^\infty$-estimates of the
kind described in the introduction. In view of our intended generalization of
the scalar case to pseudodifferential operators, let us simply focus on this
situation in the Cauchy problem (\ref{PsDE}-\ref{initialvalue}). We have
$f\in\G(\R^{n+1})$ and $g\in\G(\R^n)$ and the spatial operator $A$ is a
differential operator of the form
\[
    A = \sum\limits_{j=1}^n a_j(x,t) \d_{x_j} + b(x,t)
\]
where the coefficients $a_j$, $b$ are in $\G(\R^{n+1})$, $a_j$ real. Note that
a generalized symbol for $A$ is given by
\beq\label{PDOsymbol}
    i \sum_{j=1}^n a_{j,\eps}(x,t) \xi_j + b_\eps(x,t),
\eeq
where $a_{j,\eps}$, $b_\eps$ are any representatives of $a_j$, $b$;
$a_{j,\eps}$ may taken to be real-valued.

Sufficient conditions for existence and uniqueness of a
solution $u\in\G(\R^{n+1})$ to (\ref{PsDE}-\ref{initialvalue}) are as follows:
\begin{enumerate}
\item $a_j$, $b$ are equal to a (classical) constant for large $|x|$ (any kind
    of uniform boundedness in $x$ and $\eps$ for large $|x|$ would do; it
    ensures uniqueness and enables one to use partition of unity arguments in
    the proof)
\item $b$ as well as $D_k a_j$ are of \emph{log-type}, i.e., the asymptotic
    norm estimates (of order $0$) have bounds $O(\log(1/\eps))$ (this ensures
    existence by guaranteeing moderateness from energy estimates).
\end{enumerate}

Counter examples show that none of the two conditions can be dropped without
losing existence or uniqueness in general.

\begin{rem} It turns out that the non-uniqueness effect as
constructed in \cite[Example 1.4]{O:89} disappears in $\Gtwo$. (In mentioned
example, the constructed initial value $v(x,0) = \cl{(\chi(x + 1/\eps))}$,
$\chi\in\D(\R)$ with $\chi(0)=1$, is negligible in $\G$ but gives
$\ltwo{v_\eps(.,0)} = \ltwo{\chi} > 0$, hence is nonzero on $\Gtwo$.) As a
matter of fact, the $L^2$-energy estimates, to be discussed in the following
section, directly yield uniqueness; this holds even with coefficients that
allow for logarithmic growth as $\eps \to 0$ throughout the entire domain.
\end{rem}

The non-locality of pseudodifferential operators seems to prohibit an adaption
of the proof technique of \cite{LO:91}, where one is able to pass from
$L^2$-energy estimates to local $L^\infty$-estimates. On the other hand, when
working in $\Gtwo$, there is also the structural advantage of having good
mapping properties of pseudodifferential operators with uniform symbol
estimates on Sobolev spaces.

\section{Preparatory energy estimates}

Our proof of unique solvability of the Cauchy problem will be based on energy
estimates, with precise growth estimates of all appearing constants depending
on the regularization parameter $\eps$ as $\eps \to 0$. This in turn is solely
encoded into the generalized symbol in form of the semi-norm estimates of the
regularizing (resp.\ defining) family of symbols. Therefore, and also to make
the structure more transparent, we will first state the preparatory estimates
for smooths symbols in terms of explicit dependencies on symbol semi-norms and
insert the $\eps$-asymptotics only later on.

In order to maintain close resemblance in notation with the cases of
differential operators or decoupled systems, we shall write the symbol of $A$
in the form $i\, a(t,x,\xi)$ with $a \in \Cinf([0,T],S^1)$; in other words, we
review energy estimates for the operator \beq\label{Pia}
    P := \d_t + i\, a(t,x,D_x)
\eeq
under the hyperbolicity assumption
\beq\label{symm_hyp}
    a = a_1 + a_0 \qquad \text{ with $a_1$ real-valued, $a_0$ of order $0$},
\eeq
or equivalently, that
\[ \leqno{(\ref{symm_hyp})'} \qquad\qquad \qquad
    a(t,x,D_x) - \adjoint{a(t,x,D_x)} \quad \text{is of order }0.
\]

Besides stating the general case in the following proposition we also give
details on two special instances. These are of interest in applications and
allow for certain improvements concerning the regularity assumptions in terms
of symbol derivatives, which are required in the constants of the basic energy
estimate.
\begin{prop}\label{energy_prop} Assume that $P$ is the operator given in
(\ref{Pia}) and such that (\ref{symm_hyp}) holds. Let
$u\in\Con([0,T],H^1(\R^n)) \cap \Con^1([0,T],L^2(\R^n))$ and define $f := P u
\in \Con([0,T],L^2(\R^n))$. Then we have the energy estimate
\begin{multline}\label{energy}
    \ltwo{u(t)}^2 \leq \ltwo{u(0)}^2 + \int_0^t \ltwo{f(\tau)}^2\, d\tau + \\
    + C \, \Big( 1 + Q^0_{0,k_n',l_n'}(a_0) + Q^1_{0,k_n,l_n}(a_1)\Big)
        \int_0^t \ltwo{u(\tau)}^2\,d\tau,
\end{multline}
where the constant $C>0$ as well as $k_n'$, $l_n'$, $k_n$, $l_n$ are
independent of $u$ and can be chosen according to certain assumptions on the
symbol $a$ as follows:
\begin{description}
\item[(a) General case] We have $k_n' = l_n' = \floor{n/2} + 1$,
    $k_n = 3 (\floor{n/2} + 1)$,
    $l_n = 2(n + 2)$ and $C$ depends   only on the dimension $n$.
\item[(b) Constant for large $|x|$] If there is
    $r_0 \geq 0$ such that
  \beq\label{a_const}
    a(t,x,\xi) = h(t,\xi) \quad \text{ whenever } |x| \geq r_0,
  \eeq
    where $h$ is a symbol of order $1$ (with parameter $t$ and no $x$
    variable), then $C$ depends only on $n$, $r_0$ and the semi-norm orders
    are at most $k_n' = 0$, $l_n' = n + 1$, $k_n = 1$,
    $l_n = n + 2$.
\item[(c) Real symbol] If in addition (to any of the assumptions above) the
symbol $a$ is real-valued, so that $a_0$ is real as well
    in (\ref{symm_hyp}), then the term $Q^0_{0,k_n',l_n'}(a_0)$ can be dropped
    in (\ref{energy}).
\end{description}
\end{prop}
\begin{proof}
Using the standard decomposition of the operator $a_1(t,x,D_x)$ into self- and
skew-adjoint part, $a_1 = (a_1 + \adjoint{a_1})/2 + (a_1 - \adjoint{a_1})/2$, 
we obtain
\begin{multline} \label{basic_energy}
\diff{t}\ltwo{u(t)}^2 = 2\, \Re \inp{\d_t u(t)}{u(t)}  \\
    = 2\, \Re \inp{f(t)}{u(t)} + 2\, \Im \inp{a_1(t,x,D_x) u(t)}{u(t)}
		+ 2\, \Im \inp{a_0(t,x,D_x) u(t)}{u(t)}	 \\
    \leq \ltwo{f(t)}^2  + \ltwo{u(t)}^2 +
        \| a_1(t,x,D_x) - \adjoint{a_1(t,x,D_x)}\|\, \ltwo{u(t)}^2 \\
	+ 2 \| a_0(x,D)\|\, \ltwo{u(t)}^2 .
\end{multline}
The operator norms (taken with respect to $L^2$) are finite by
(\ref{symm_hyp}) and (\ref{symm_hyp})$'$ (cf.\ \cite[Theorem
18.1.11]{Hoermander:V3} or \cite[Ch.2, Theorem 4.1]{Kumano-go:81}) and we will
derive explicit estimates for these. 

For the proof of \emph{case (c)}, observe that we have 
$\Im \inp{a_0(t,x,D_x) u(t)}{u(t)} = 0$ if $a_0$ is real, so that the last
term on te right-hand side of (\ref{basic_energy}) can be dropped from all
further considerations.

\emph{Case (a):} We use a representation of the symbol of $b_1(t,x,D_x) :=
a_1(t,x,D_x) - \adjoint{a_1(t,x,D_x)}$ with integral remainder
terms as it is 
developed in \cite[Ch.2,1-3]{Kumano-go:81} or \cite[Ch.1,5-6]{Cordes:95}.
According to this 
(or as sketched in the Appendix below),
the zero order symbol $b_1(t,x,\xi)$ is given by
 \begin{multline}\label{expansion}
    b_1(t,x,\xi) :=  a_1(t,x,\xi) - \ovl{a_1(t,x,\xi)}
        - i \sum_{j = 1}^n \int_0^1 r_{j,\theta}(t,x,\xi)\, d\theta \\
	=  - i \sum_{j = 1}^n \int_0^1 r_{j,\theta}(t,x,\xi)\, d\theta,
 \end{multline}
since $a_1(t,x,\xi)$ is real-valued, where
 \beq\label{remainder}
    r_{j,\theta}(t,x,\xi) = \iint e^{-i y\cdot\eta} \,
        \ovl{\d_{\xi_j} \d_{x_j} a_1(t,x+y,\xi + \theta\eta)}\, dy \dq\eta
 \eeq
 in the sense of oscillatory integrals. 

 As a close inspection of the proof of \cite[Ch.2, Lemma 2.4]{Kumano-go:81}
 shows (which we detail in the Appendix), we have the following estimate for
 all 
 $d\in\N_0$, $(\al,\be)\in\N_0^{2n}$
 \beq\label{remainder_estimate}
    |\d_t^d \d_\xi^\al \d_x^\be r_{j,\theta}(t,x,\xi)|
        \leq C_{d,\al,\be}\, (1 + |\xi|)^{-|\al|}\,
            Q^1_{d, n + 2 + |\al|, n + 2 + |\al| + |\be|}(a_1),
 \eeq
which is uniform with respect to $\theta\in[0,1]$. Combined with formula
(\ref{expansion}) for $b_1$ this yields
\[
    Q^0_{0,k,l}(b_1 + 2 a_0) \leq 2\, Q^0_{0,k,l}(a_0) + C_{k,l}\,
        Q^1_{0, n + 2 + k, n + 2 + k + l}(a_1).
\]
By the theorem of Calder\'{o}n-Vaillancourt (or one of its variants, cf.\
\cite[Ch.I, Th\'{e}or\`{e}me 3]{CM:78}, \cite{Hounie:86}, \cite[Ch.3, Corollary
1.3]{Cordes:95}), we have the general $L^2$-operator norm estimate
\[
    \| b(t,x,D_x) \| \leq C'_n Q^0_{0, \floor{n/2} + 1,
    \floor{n/2} + 1}(b) 
\]
whenver $b\in\Cinf([0,T],S^0)$. 
Therefore we conclude that
\[
    \| b_1(t,x,D_x)\| + 2 \| a_0(x,D) \| 
	\leq C_n \Big( Q^0_{0,\floor{n/2}+1,\floor{n/2}+1}(a_0) 
		+ Q^1_{0,k,l}(a_1) \Big)
\]
for any $k \geq n + \floor{n/2} + 3$ and $l \geq n + 2 \floor{n/2} + 4$. This
completes the proof of the general case.

\emph{Case (b):} Let $\chi\in\D(\R^n)$ such that $\chi(x) = 1$ for $|x|\leq
r_0$ and $0 \leq \chi  \leq 1$. Then the term $a(t,x,\xi) - \ovl{a(t,x,\xi)}$
occurring in (\ref{expansion}) can be written in the form
\begin{multline*}
    \chi(x) (a_0(t,x,\xi) - \ovl{a_0(t,x,\xi)}) +
        (1 - \chi(x)) (h(t,\xi) - \ovl{h(t,\xi)}) \\
            := b_0(t,x,\xi) + (1-\chi(x)) h_0(t,\xi).
\end{multline*}
The second part in this decomposition is the operator symbol of a convolution
with bounded symbol (since $h_0 := h - \ovl{h}$ is of order $0$), composed with
multiplication by $1-\chi$ from the left. Hence the $L^2$ operator norm
corresponding to this second summand has the following upper bound
\[
    \| (1-\chi(x)) h_0(t,D_x) \| \leq \linf{1-\chi} \linf{h_0} \leq 2 \,\linf{a_0}.
\]
Note that $b_0$ is a symbol of order zero with support contained in $|x| \leq
r_0$. We will estimate the operator norm of $B_0 := b_0(t,x,D_x)$ (on
$L^2(\R^n_x)$, uniformly with respect to $t\in[0,T]$) via the Schwartz kernel
$\wt{K_0}$ of the ``Fourier transformed'' operator $\wt{B_0} := \F\circ B_0
\circ\F\inv$ and using the fact that
 \[
    \| B_0 \| = (2\pi)^{n/2} \| \wt{B_0} \|.
 \]
As a distribution in $\Cinf([0,T],\S'(\R^{2n}))$, the kernel
 is computed from the symbol by the formula
 \beq\label{FT_Kernel}
    \wt{K_0}(t,\xi,\eta) = (2\pi)^{-n}\, \F\big(b_0(t,.,\eta)\big)(\xi-\eta).
 \eeq
Since $x\mapsto b_0(t,.,\eta)$ has compact support it follows that $\wt{K_0}$
is smooth on $[0,T]\times\R^{2n}$; in fact, we will see that it is an
integrable kernel and hence we may apply a classical lemma of Schur (cf.\
\cite[Lemma 18.1.12]{Hoermander:V3}). Before doing so, we will first show that
the remainder terms in (\ref{expansion}) are of a similar form.

Consider formula (\ref{remainder}) and introduce the short-hand notation $b_j
:= \d_{\xi_j} \d_{x_j} \ovl{a}$. Then $z\mapsto b_j(t,z,\zeta)$ has compact
support in $|z| \leq r_0$ and we may write
 \[
    r_{j,\theta}(t,x,\xi) =
        \F\inv_{\eta\to x}\Big(\F\big(b_j(t,.,\xi+\theta \eta)\big)(\eta)\Big).
 \]
Now let $R_{j,\theta} := r_{j,\theta}(t,x,D_x)$ and define, exactly as above,
the corresponding operator $\wt{R_{j,\theta}}$ with intertwining Fourier
transforms; denote by $\wt{K_{j,\theta}}$ its Schwartz kernel. The above
representation for the symbol $r_{j,\theta}$ in terms of $b_j$ and direct
computation yields the formula
 \beq\label{Kj_formula}
    \wt{K_{j,\theta}}(t,\xi,\eta) = (2\pi)^{-n}\,
        \F\big(b_j(t,.,\eta + \theta(\xi - \eta))\big)(\xi-\eta).
 \eeq

Equations (\ref{FT_Kernel}) and (\ref{Kj_formula}) have the following structure
in common: we have a symbol $d\in\Cinf([0,T],S^0)$ which vanishes when $|x|\geq
r_0$ and a smooth kernel $\wt{K}$ defined by
\[
    \wt{K}(t,\xi,\eta) := (2\pi)^{-n}\,
    \F\big(d(t,.,f(\xi,\eta))\big)(\xi-\eta),
\]
where $f \colon \R^{2n} \to \R^n$ is a linear map. In order to apply Schur's
lemma we estimate the partial $L^1$-norms of the kernel and obtain
\begin{multline*}
    (2\pi)^n \int\limits_{\R^n} |\wt{K}(t,\xi,\eta)| \, d\xi =
        \int\limits_{\R^n} |\F\big(d(t,.,f(\xi,\eta))\big)(\xi-\eta)|\, d\xi \\
        = \int\limits_{\R^n} |\F\big(d(t,.,f(\xi+\eta,\eta))\big)(\xi)|\,d\xi
        \leq
        \int\limits_{\R^n}\sup_{t,\zeta}|\F\big(d(t,.,\zeta)\big)(\xi)|\,d\xi,
\end{multline*}
and similarly
\[
    (2\pi)^n \int\limits_{\R^n} |\wt{K}(t,\xi,\eta)| \, d\eta \leq
    \int\limits_{\R^n}\sup_{t,\zeta}|\F\big(d(t,.,\zeta)\big)(\eta)|\,d\eta.
\]

\emph{Assertion:} There exists a constant $c(n,r_0)$, depending only on $n$ and
$r_0$, such that
 \beq\label{int_estimate}
    \int\limits_{\R^n} \sup_{t,\zeta} | \F\big(d(t,.,\zeta)\big)(\mu) | \, d\mu
        \leq c(n,r_0) \, Q^0_{0,n+1,0}(d).
 \eeq
Noting that $\int |\d_x^\be d(t,x,\zeta)|\, dx \leq c_n r_0^n 
\linf{\d_x^\be d(t,.,\zeta)}$, 
the proof is exactly as in \cite[Theorem 18.1.11${}'$]{Hoermander:V3}. 

In summary, applying (\ref{int_estimate}) and the general integral kernel
estimates above to the kernels given by (\ref{FT_Kernel}) and
(\ref{Kj_formula}) (note that $b_j$ involves first-order derivatives of $a$ in
$x$ and $\xi$ already) we have proved the claims of case (b) in the
proposition.
\end{proof}

\section{Colombeau solutions}

We return to the scalar pseudodifferential Cauchy problem
\begin{align}
    \d_t u + A u &=  f  \quad \text{in $X_T$}, \label{C_PsDE}\\
    u(0) &= g, \label{C_initialvalue}
\end{align}
where $X_T = \R^n\times (0,T)$ and with data $f\in\Gtwo(X_T)$ and
$g\in\Gtwo(\R^n)$. $A$ is a generalized pseudodifferential operator of order
$1$. More precisely, we assume that
 \[
    A \colon \Gtwo(X_T) \to \Gtwo(X_T) \quad \text{is given by}\quad
     (u_\eps)_\eps \mapsto (i\, a_\eps(t,x,D_x)u_\eps)_\eps,
 \]
where $(a_\eps(t,x,\xi))_\eps$ is a generalized symbol of order $1$ with
parameter $t\in[0,T]$. In addition, we impose the hyperbolicity assumption
 \beq \label{C_symm_hyp}
    \forall \eps\in(0,1] \colon\quad a_\eps = a_{1,\eps} + a_{0,\eps}
    \qquad \text{ $a_{1,\eps}$ real-valued, $a_{0,\eps}$ of order $0$}.
 \eeq

The semi-norms in the basic energy estimate (\ref{energy}) now depend on
$\eps\in (0,1]$, and, upon applying Gronwall's inequality, will appear as
exponents in the $L^2$-norm estimates of a prospective generalized solution;
this suggests to assume logarithmic bounds on the symbols. We say that a
generalized symbol $(b_\eps)_\eps$ of order $m$ (with parameter $t\in[0,T]$) is
of \emph{log-type up to order (k,l)} if
 \beq\label{log_type}
    Q^m_{0,k,l}(b_\eps) = O(\log(1/\eps)) \qquad (\eps\to 0).
 \eeq

\begin{thm}\label{main_thm} Let $A$ be a generalized first-order
pseudodifferential operator, defined by the generalized symbol
$(ia_\eps)_{\eps\in(0,1]}$ with parameter $t\in[0,T]$, and satisfying the
hyperbolicity assumption (\ref{C_symm_hyp}). Assume that $(a_{1,\eps})_\eps$ 
is of log-type up to order $(k_n,l_n + 1)$ and that $(a_{0,\eps})_\eps$ is of
log-type up to order $(k_n',l_n')$.

Then for any given $f\in\Gtwo(X_T)$, $g\in\Gtwo(\R^n)$ the Cauchy problem
(\ref{C_PsDE}-\ref{C_initialvalue}) has a unique solution $u\in\Gtwo(X_T)$ 
if $k_n = 3(\floor{n/2} + 1)$, $l_n = 2(n + 2)$, $k_n' = l_n' = \floor{n/2} +
1$

Furthermore, we have variants of the log-type requirements in the following
two cases:
 \begin{enumerate}
 \item If there is $r_0 \geq 0$ and an $x$-independent generalized symbol
     $(h_\eps(t,\xi))_\eps$ such that
     \beq\label{C_a_const}
        a_\eps(t,x,\xi) = h_\eps(t,\xi) \qquad \text{ when $|x|\geq r_0$},
     \eeq
     then we may put $k_n = 1$, $l_n = n + 2$, $k_n' = 0$,
     $l_n' = n +1 $.
 \item If $a_\eps$ is real-valued for every $\eps\in(0,1]$ then no
    log-type assumption on $a_{0,\eps}$ is required.
 \end{enumerate}
\end{thm}
\begin{proof}
Let $(g_\eps)_\eps\in g$, $(f_\eps)_\eps\in f$ be representatives. At fixed,
but arbitrary, $\eps\in(0,1]$ we consider the smooth Cauchy problem
\begin{align}
    \d_t u_\eps + i a_\eps(t,x,D_x) u_\eps &=  f_\eps  \quad \text{in $X_T$},
        \label{eps_PsDE}\\
    u_\eps(0) &= g_\eps. \label{eps_initialvalue}
\end{align}
It has a unique solution $u_\eps\in\Cinf([0,T],H^\infty(\R^n))$, thus
constituting a solution candidate $(u_\eps)_\eps$ (cf.\ \cite[Ch.7, Theorem
3.2]{Kumano-go:81} or \cite[Ch.6, Theorem 2.1]{Cordes:95} with additional
$t$-regularity following directly from the equation). We have to show that
$(u_\eps)_\eps \in \E_{M,L^2}(X_T)$.

Denote by $C_\eps := C( 1 + Q^0_{0,k_n',l_n'}(a_{0,\eps}) +
Q^1_{0,k_n,l_n}(a_\eps))$ the constant occurring in the energy estimate
(\ref{energy}) applied to $u_\eps$. Gronwall's lemma implies
 \beq\label{Gronwall_energy}
    \ltwo{u_\eps(t)}^2 \leq
        \Big( \ltwo{g_\eps}^2 + \int_0^T \ltwo{f_\eps(\tau)}^2\, d\tau \Big)\,
        \exp(C_\eps T).
 \eeq
By hypothesis we have $C_\eps = O(\log(1/\eps))$ as $\eps\to 0$. Thus we obtain
\emph{uniqueness} immediately from (\ref{Gronwall_energy}) -- once moderateness
is established -- because null family estimates for $f_\eps$, $g_\eps$ then
imply such for $u_\eps$ as well.

For the proof of \emph{existence}, we first observe that the basic estimate for
$\norm{u_\eps}{L^2(X_T)} \leq T \sup_{t\in[0,T]} \ltwo{u_\eps(t)} =
O(\eps^{-N})$ follows at once from (\ref{Gronwall_energy}) by the moderateness
of the data. It remains to prove moderateness estimates for the higher order
derivatives of $u_\eps$.

\emph{$x$-derivatives:} Let $0 \not= \al\in\N_0^n$ and apply $\d_x ^\al$ to
equation (\ref{eps_PsDE}). It follows by induction and simple commutator
relations of $a_\eps(t,x,D_x)$ with $\d_{x_j}$ that this produces an equation
of the following structure. Denote by $e_j = (\de_{j,k})_{k=1}^n$ the
$j^{\text{th}}$ standard basis vector in $\R^n$ then
 \beq\label{dxalpha}
    \d_t\d_x^\al u_\eps + i a_\eps(t,x,D_x) \d_x^\al u_\eps +
        i \sum_{1 \leq j \leq n \atop \al_j \not= 0} (\d_{x_j}a_{1,\eps})(t,x,D_x)
        \d_x^{\al-e_j} u_\eps = F_{\eps,\al},
 \eeq
 where $F_{\eps,\al}$ equals the sum of $\d_x^\al f_\eps$ plus, if $ |\al| \geq 2$,
 a linear combination of terms of the form
  \beq\label{Falpha2}
    (\d_x^\be a_\eps)(t,x,D_x) \d_x^{\al-\be} u_\eps \qquad
        \text{ with $\be \leq \al$ and $2 \leq |\be|$},
  \eeq
 and
 \beq\label{Falpha3}
    (\d_{x_j}a_{0,\eps})(t,x,D_x)\d_x^{\al-e_j} u_\eps \qquad
        \text{ where } \al_j \not= 0.
 \eeq
Assume that moderateness of $\ltwo{\d_x^\ga u_\eps}$ has been established
already when $|\ga| < |\al|$. Since $\d_x^\be a_\eps$ is of order $1$ we have
 \begin{multline*}
   \ltwo{(\d_x^\be a_\eps)(t,x,D_x) \d_x^{\al-\be} u_\eps(t)} \leq C_1\,
    Q^1_{0,m,m}(\d_x^\be a_\eps)\, \norm{\d_x^{\al-\be} u_\eps(t)}{H^1} \\
   \leq C_1' \, Q^1_{0,m,m}(\d_x^\be a_\eps)\,
    \max_{|\ga| < |\al|} \ltwo{\d_x^\ga u_\eps(t)},
 \end{multline*}
where $C_1$, $C_1'$, and $m$ depend only on the dimension $n$ (\cite[Ch.3,
Theorem 2.7]{Kumano-go:81}). Similarly, since $\d_{x_j} a_{0,\eps}$ is of order
$0$ we also have
 \begin{multline*}
   \ltwo{(\d_{x_j} a_{0,\eps})(t,x,D_x) \d_x^{\al-e_j} u_\eps(t)} \leq C_2\,
    Q^0_{0,m',m'}(\d_{x_j} a_{0,\eps})\, \ltwo{\d_x^{\al-e_j} u_\eps(t)} \\
   \leq C_2' \, Q^0_{0,m',m'}(\d_{x_j} a_{0,\eps})\,
    \max_{|\ga| < |\al|} \ltwo{\d_x^\ga u_\eps(t)},
 \end{multline*}
where $C_2$, $C_2'$, depend only on the dimension and $m' = \floor{n/2}+1$.
 Hence, by the induction
hypothesis, we have $\ltwo{F_{\eps,\al}(t)} = O(\eps^{-N})$ as $\eps\to 0$
uniformly in $t\in[0,T]$ for some $N$.

We return to equation (\ref{dxalpha}), consider it as an equation for $v_\eps
:= \d_x^\al u_\eps$, and supply the initial value $v_\eps(0) = \d_x^\al
u_\eps(0) = \d_x^\al g_\eps$. Applying the basic technique from the beginning
of the proof of the energy estimate (\ref{energy}) to equation (\ref{dxalpha})
we obtain
 \begin{multline*}
    \diff{t} \ltwo{v_\eps(t)}^2 = 2 \, \Re \inp{\d_t v_\eps(t)}{v_\eps(t)}
    \\ \leq \ltwo{F_{\eps,\al}(t)}^2 +
      \ltwo{v_\eps(t)}^2 +
      \| a_\eps(t,x,D_x) - \adjoint{a_\eps(t,x,D_x)}\| \ltwo{v_\eps(t)}^2 \\
    + \sum_{1\leq j\leq n \atop \al_j \not= 0}
    \|(\d_{x_j}a_{1,\eps})(t,x,D_x) - \adjoint{(\d_{x_j}a_{1,\eps})(t,x,D_x)}\|
    \big( \ltwo{\d_x^{\al - e_j}u_\eps(t)}^2 +  \ltwo{v_\eps(t)}^2 \big).
 \end{multline*}
If we define $G(\d_x a_{1,\eps})(t) := \sum_{j=1}^n
    \|(\d_{x_j}a_{1,\eps})(t,x,D_x) - \adjoint{(\d_{x_j}a_{1,\eps})(t,x,D_x)}\|$
then we get
 \begin{multline*}
    \diff{t} \ltwo{v_\eps(t)}^2 \leq \ltwo{F_{\eps,\al}(t)}^2 +
        G(\d_x a_{1,\eps})(t) \max_{|\ga| < |\al|} \ltwo{\d_x^\ga u_\eps(t)}^2 \\
    + \Big( 1 + \| a_\eps(t,x,D_x) - \adjoint{a_\eps(t,x,D_x)}\|
    + G(\d_x a_{1,\eps})(t) \Big) \ltwo{v_\eps(t)}^2.
 \end{multline*}
The term $H_{\eps,\al}(t) := \ltwo{F_{\eps,\al}(t)}^2 + G(\d_x
a_{1,\eps})(t)\cdot \max_{|\ga| < |\al|} \ltwo{\d_x^\ga u_\eps(t)}^2$ is of
moderate growth and, by the proof of Proposition \ref{energy_prop}, we have
 \begin{multline*}
    1 + \| a_\eps(t,x,D_x) - \adjoint{a_\eps(t,x,D_x)}\| + G(\d_x
    a_{1,\eps})(t)\\
    \leq  C\, \big( 1 + Q^0_{0,k_n',l_n'}(a_{0,\eps}) + Q^1_{0,k_n,l_n}(a_\eps)
        + \sum_{j=1}^n Q^1_{0,k_n,l_n}(\d_{x_j} a_{1,\eps}) \big) :=
        \wt{C_\eps},
 \end{multline*}
 which is a log-type constant by the hypotheses of the theorem. Note that the
specifications of $k_n$, $l_n$, $k_n'$, $l_n'$ for the general case match those
of Proposition \ref{energy_prop}, case (a), whereas the hypotheses in (i), (ii)
match cases (b), (c) there.  Thus, we prove all assertions of the theorem
simultaneously when the notation is understood in this way. Finally,
integration with respect to $t$ and Gronwall's lemma yield the estimate
  \[
    \ltwo{\d_x^\al u_\eps(t)}^2 \leq \Big( \ltwo{\d_x^\al g}^2 +
    \int_0^T H_{\eps,\al}(\tau)\, d\tau \Big) \exp(\wt{C_\eps} T) = O(\eps^{-M})
  \]
 for some $M$ and $\eps$ sufficiently small. Hence $\norm{\d_x^\al
 u_\eps}{L^2(X_T)}$ satisfies a similar estimate. In particular, we have the
 same bounds on the spatial Sobolev norms $\norm{u_\eps(t)}{H^k}$ for $k$
 arbitrary and uniformly in $t\in[0,T]$.

\emph{$t$- and mixed derivatives:} Equations (\ref{eps_PsDE}) and
(\ref{dxalpha}) directly imply estimates of the form $\ltwo{\d_t \d_x^\al
u_\eps(t)} = O(\eps^{-N})$ for any $\al\in\N_0^n$ (and uniformly in $t$). To
proceed to higher order $t$-derivatives, we simply differentiate equations
(\ref{eps_PsDE}), resp.\ (\ref{dxalpha}), with respect to $t$. The Sobolev
mapping properties of the operators $(\d_t \d_x^\be) a_\eps(t,x,D_x)$ and
moderateness assumptions on the symbols then yield the desired estimates for
$\ltwo{\d_t^l \d_x^\al u_\eps(t)}$ successively for $l=0,1,2,\ldots$ and $\al$
arbitrary.
\end{proof}

\begin{rem}
(i) The key assumptions in Theorem \ref{main_thm} are the log-type estimates on
the symbol. We know already from the differential operator case that they
cannot be dropped completely. However, these are sufficient conditions and
merely reflect the various operator norm bounds available for zero order
symbols (as used in proving the energy estimates). Thus they cannot be expected
to be sharp. In fact, the value of the theorem lies in a general feasibility
proof and any special structure inherent in a concrete symbol under
consideration in applications might allow for improvement.

(ii) In order to meet the log-type conditions of the above theorem in a
specific symbol regularization one may call on a re-scaled mollification as
described in \cite{O:89}. To illustrate this procedure, let us assume that the
non-smooth symbol of order $m$ is given as the measurable bounded function
$a(x,\xi)$ such that for almost all $x$ the partial function $\xi \mapsto
a(x,\xi)$ is smooth and satisfies for all $\al\in\N_0^\al$ an estimate
 \[
    \linf{\d_\xi^\al a(.,\xi)} \leq C_\al (1 + |\xi|)^{m - |\al|}.
 \]
\end{rem}
Let $\rho$ be a mollifier and let $0 < \om_\eps \leq (\log(1/\eps))^{1/k}$ for
some $k\in\N$, $\om_\eps \to \infty$ as $\eps\to 0$. Let $\rho^\eps(y) :=
\om_\eps^n \rho(\om_\eps y)$ and define the regularized symbol by
$a_\eps(x,\xi) := (\rho^\eps * a_\eps(.,\xi))(x)$ (convolution with respect to
the $x$-variable only). Then it is easy to check that $a_\eps\in S^m$ and of
log-type up to order $(\infty,k)$.

As in \cite{LO:91}, essentially by inspection of the above existence proof,
we establish compatibility with distributional or smooth solutions, that is
macroscopic regularity in a certain sense, when the symbol is smooth. 
\begin{cor}\label{compatibility_cor} In Theorem \ref{main_thm}, assume that
$A$ is given by a smooth symbol $a\in\Cinf([0,T],S^1)$, i.e., $a_\eps = a$ for
all $\eps\in(0,1]$.  
\begin{enumerate}
\item If $f\in H^\infty(X_T)$, $g\in H^\infty(\R^n)$, then the generalized
solution $u\in\Gtwo(X_T)$ is equal to the classical smooth solution.
\item Let $f\in C([0,T],H^s(\R^n))$ and $g\in H^s(\R^n)$ for some $s\in\R$,
and $v$ be the unique distributional solution to
(\ref{PsDE}-\ref{initialvalue}) in $C([0,T],H^s)$. Define generalized data
for problem (\ref{C_PsDE}-\ref{C_initialvalue}) by  
$\wt{f} := \cl{(f_\eps)_\eps}\in\Gtwo(X_T)$ (resp.\
$\wt{g} := \cl{(g_\eps)_\eps}\in\Gtwo(\R^n)$), where $f_\eps\in H^\infty(X_T)$
(resp.\ $g_\eps\in H^\infty(\R^n)$) are moderate regularizations such that
$f_\eps \to f$ in $C([0,T],H^s(\R^n))$ (resp.\ $g_\eps \to g$ in $H^s(\R^n)$)
as $\eps \to 0$. 
If $u$ is the corresponding generalized solution in $\Gtwo(X_T)$ then it is
associated with the distributional solution $v$. 
\end{enumerate}
\end{cor}
\begin{proof}
\emph{Part (i):} Since we may choose the constant nets $(f)_\eps$,
$(g)_\eps$ as representatives of the classes of $f$ and $g$ in $\Gtwo$, and
$a_\eps = a$ by assumption, we obtain the classical smooth solution to
equation (\ref{eps_PsDE}-\ref{eps_initialvalue}) as a representative of the
unique Colombeau solution.

\emph{Part (ii):} The unique solution $v\in C([0,T],H^s)$ to
(\ref{PsDE}-\ref{initialvalue}) depends continuously on the data 
$f$ and $g$ by the closed graph theorem. Hence the solution representative
$u_\eps$, defined as the solution to (\ref{eps_PsDE}-\ref{eps_initialvalue}), 
converges to $v$ in $C([0,T],H^s)$ as $\eps \to 0$.
\end{proof}

Finally, we prove that the intrinsic regularity property for the generalized
solution holds if the data are in $\Gtwo^\infty$ and the generalized symbol is
only mildly generalized, namely satisfies additional slow scale conditions. 
This notion was introduced and investigated in some detail in \cite{HO:02} and
found to be crucial for regularity theory of partial differential equations. 
Recall that a net $(r_\eps)_\eps$ of complex numbers is said to be of
\emph{slow scale} if it satisfies   
 \[
   \forall p\geq 0 \colon \,\, |r_\eps|^p = O(1/\eps) \quad (\eps \to 0).
 \]
In the proposition below, we call a net $(s_\eps)_\eps$ of complex numbers a
\emph{slow-scale log-type net} if there is a slow scale net $(r_\eps)_\eps$ of
real numbers, $r_\eps \geq 1$, such that 
 \[
	|s_\eps| = O(\log(r_\eps)) \quad (\eps \to 0).
 \]
\begin{prop}\label{reg_prop} In Theorem \ref{main_thm}, assume all
log-type conditions to be replaced by slow-scale log-type estimates and, in
addition, that $(a_\eps)_\eps$ is of slow scale in each derivative. By the
latter, we mean that for all $j$, $k$, $l$, we can find a slow scale net 
$(r_\eps)_\eps$ positive real numbers such that 
 \[ 
 	Q^1_{j,k,l}(a_\eps) = O(r_\eps) \qquad (\eps \to 0). 
 \] 
Then $f\in\Gtwo^\infty(X_T)$ and $g\in\Gtwo^\infty(\R^n)$ implies
$u\in\Gtwo^\infty(X_T)$. In particular, this is always true when the symbol of
$A$ is smooth (as in the Corollary above).
\end{prop}
\begin{proof} Thanks to the explicit assumptions this is straightforward 
by an inspection of the proof of Theorem \ref{main_thm}. To be more precise,
assume that we have a uniform $\eps$-growth, say $\eps^{-M}$, for the
derivatives of $f$ and $g$; i.e., for 
all $k$, $\al$, we have
 \[ 
 	\ltwo{\d_t^k \d_x^\al f_\eps} = O(\eps^{-M}) \quad\text{ as well as }\quad
	\ltwo{\d_x^\al g_\eps} = O(\eps^{-M}).
 \]
Note that all constants involving $a_\eps$ in the energy and Sobolev estimates
throughout the proof yield only slow scale factors. (Observe again, that in the
exponential factors in all energy estimates we only need a fixed finite order
of derivatives, corresponding to $k_n'$, $l_n'$ etc.) Thus, the same induction
argument shows that we obtain for all $k$, $\al$ a certain slow scale net
$(r_\eps)_\eps$ of positive real numbers, such that
 \[
	\ltwo{\d_t^k \d_x^\al u_\eps} = O(r_\eps \eps^{-M}) = O(\eps^{-M-1}),
 \] 
whiches proves the assertion. 
\end{proof}

\begin{rem} 
(i) The somewhat extensive slow-scale log-type conditions in the above
proposition are by far not necessary for regularity, but are suited to make
the energy estimates, with their exponential constants, directly applicable.
We expect that these can be relaxed at least to plain slow scale conditions by
appealing to pseudodifferential parametrix techniques (cf.\
\cite{Garetto:02,GGO:03}).

(ii) A slow-scale property of $(a_\eps)_\eps$ (in all derivatives) is
implied, for example,  
by the log-type assumptions on $(a_\eps)_\eps$ if, in addition, only
a $\Ginf$-type regularity of $(a_\eps)_\eps$ is assumed. This follows from
\cite[Proposition 1.6]{HOP:03} and the fact that $\log(1/\eps)$ is a slow scale
net. 
\end{rem}

\appendix{{\bf Appendix: Remainder term estimates}}

We briefly outline a proof of (\ref{expansion}) and verify the precise form of
the estimate (\ref{remainder_estimate}); it is an adaption
of the reasoning in \cite[Ch.2, Sections
2-3]{Kumano-go:81}; thereby, we also recall the precise meaning of the
oscillatory integral (\ref{remainder}). We may suppress the dependence of all
symbols on the parameter $t$, since it will be clear that all steps in the
process respect continuity (or smoothness) with respect to it and yield uniform
bounds in all estimates when $t$ varies in $[0,T]$.

Let $a(x,\xi)$ be a (smooth) symbol of order $1$. The starting point is the
following formula for the adjoint of $a(x,D)$, e.g. valid for $u\in\S(\R^n)$ as
iterated integral,
 \[
    \adjoint{a(x,D)}u(x) = \iint e^{i(x-y)\eta} \ovl{a(y,\eta)} u(y) \, dy
    \dq\eta.
 \]
Writing $u(y) = \int e^{iy\xi}\, \FT{u}(\xi)\, \dq\xi$ we obtain, now in the
sense of oscillatory integrals,
 \begin{multline*}
    \adjoint{a(x,D)}u(x) = \int e^{ix\xi}\, \FT{u}(\xi)
        \iint e^{i(x-y)\eta} \ovl{a(y,\xi+\eta)} \, dy\dq\eta
        \dq\xi \\
    =: \int e^{ix\xi} \, \FT{u}(\xi) a^*(x,\xi)\, \dq\xi,
 \end{multline*}
by which we define the symbol $a^*$. Using Taylor expansion
 \[
    \ovl{a(y,\xi+\eta)} = \ovl{a(y,\xi)} +
        \int_0^1 \eta \cdot \d_\xi \ovl{a(y,\xi + \theta\eta)} \, \theta
 \]
and (the oscillatory integral interpretation of) the Fourier identity $\F\inv(
\F(\ovl{a(.,\xi)}))(x) = \ovl{a(x,\xi)}$ leads to
 \[
    a^*(x,\xi) = \ovl{a(x,\xi)} + \int_0^1
        \iint \eta \cdot e^{i(x-y)\eta} \ovl{\d_\xi a(y,\xi + \theta\eta)}
        \,dy\dq\eta
    d\theta.
 \]
Noting that $\eta e^{i(x-y)\eta} = D_y (e^{i(x-y)\eta})$ and integrating by
parts yields equations (\ref{expansion}) and (\ref{remainder}). We use the
notation $\d_x \d_\xi = \sum_{j=1}^n \d_{x_j} \d_{\xi_j}$, $r_\theta = \sum_j
r_{j,\theta}$, and recall that (\ref{remainder}) can be defined as the
classical integral
 \[
     r_\theta(x,\xi) = \iint e^{-iy\eta} (1+|y|^2)^{-\la} (1-\Delta_\eta)^\la
        \big( \d_x \d_\xi \ovl{a(x+y,\xi+\theta\eta)} \big)
        \,dy\dq\eta,
 \]
 where $\la > n/2$, so that $s_\theta(x,\xi;y,\eta) := (1+|y|^2)^{-\la}
 (1-\Delta_\eta)^\la\big( \d_x \d_\xi \ovl{a(x+y,\xi+\theta\eta)} \big)$
is integrable. We have the estimate
 \begin{multline}\label{stheta_estimate}
    |\d_\xi^\al \d_x^\be s_\theta(x,\xi;y,\eta)| \leq \\
        c_{n,\al,\be} \, q^1_{2\la + 1 + |\al|,1 + |\be|}(a)\,
        (1+|y|)^{-2\la}
        (1 + |\xi + \theta\eta|^2)^{- |\al|/2},
 \end{multline}
where $c_{n,\al,\be}$ is uniform in $\theta\in[0,1]$. In order to prove
(\ref{remainder_estimate}) we have to estimate
\begin{multline*}
    \d_\xi^\al \d_x^\be r_\theta(x,\xi) = \iint e^{-iy\eta}
        \d_\xi^\al \d_x^\be s_\theta(x,\xi;y,\eta) \,dy\dq\eta \\
        = \iint\limits_{ |\eta| \leq |\xi|/2 } \!\!\!\ldots \;\;\; +
            \iint\limits_{|\eta| \geq |\xi|/2 }
            \!\!\!\ldots\;\;\;
            =: I_1 + I_2.
\end{multline*}
For an upper bound of $I_1$ we use (\ref{stheta_estimate}) and the implication
$|\eta| \leq |\xi|/2 \Rightarrow |\xi+\theta\eta| \geq |\xi|/2$ (when
$\theta\in[0,1]$) to find
 \[
    |I_1| \leq c_{n,\al,\be,\la} \,q^1_{2\la + 1 + |\al|,1 + |\be|}(a)\,
        (1 + |\xi|)^{-|\al|},
 \]
 uniformly in $\theta$. To estimate $I_2$, we first use that
 $e^{-iy\eta} = |\eta|^{-2l} (-\Delta_y)^l(e^{-iy\eta})$ and integrate by parts
 to obtain
  \[
    |I_2| \leq \int\limits_{\R^n} \int\limits_{|\eta|\geq|\xi|/2}
        |\eta|^{-2l} |(-\Delta_y)^l \d_\xi^\al \d_x^\be s_\theta(x,\xi;y,\eta)|
        \, dy\dq\eta.
  \]
We apply (\ref{stheta_estimate}) with $\be$ replaced by $\be + 2l e_j$
($j=1,\ldots,n$) to the integrand and arrive at
 \[
    |I_2| \leq c_{n,\al,\be,l,\la} \, q^1_{2\la + 1 + |\al|,2l + 1 + |\be|}(a)
        \int\limits_{|\eta|\geq|\xi|/2} |\eta|^{-2l} (1 + |\xi
        +\theta\eta|^2)^{-|\al|/2 } \, \dq\eta.
 \]
By Peetre's inequality $(1+|\xi+\theta\eta|^2)^{-|\al|/2} \leq 2^{|\al|/2}
(1+|\xi|^2)^{-|\al|/2} (1 + |\theta\eta|^2)^{|\al|/2}$, so that
 \[
    |I_2| \leq c'_{n,\al,\be,l,\la} \, q^1_{2\la + 1 + |\al|,2l + 1 + |\be|}(a)
        (1 + |\xi|)^{-|\al|} \int\limits_{|\eta|\geq|\xi|/2} |\eta|^{-2l}
        (1 + |\theta\eta|^2)^{|\al|/2} \, \dq\eta,
 \]
 where the remaining integral is finite if $2l > n + |\al|$.

 Summing up, and combining the conditions $2\la > n$, $2l > n + |\al|$, we
 have shown (\ref{remainder_estimate}).

\bibliographystyle{abbrv}
\bibliography{gueMO}

\end{document}